\documentclass{amsart}
\usepackage{cite}
\usepackage{mathpazo}
\usepackage[english]{babel}
\usepackage{amssymb}
\usepackage{latexsym}
\usepackage{xypic}
\usepackage[all]{xy}
\usepackage{multirow}
\usepackage[utf8]{inputenc}
\usepackage{pgfplots}
\usepackage{float}
\usepackage{tikz}
\usepackage{ mathrsfs }
\usepackage{tikz-3dplot}
\usepackage{subcaption}
\usetikzlibrary{
  knots,
  hobby,
  decorations.pathreplacing,
  shapes.geometric,
  calc}

\usepackage[colorinlistoftodos,prependcaption,textsize=tiny]{todonotes}
\newtheorem{theorem}{Theorem}[section]
\newtheorem{lemma}[theorem]{Lemma}

\theoremstyle{definition}
\newtheorem{definition}[theorem]{Definition}

\theoremstyle{remark}
\newtheorem{remark}[theorem]{Remark}

\newcommand{\Es}{\mathbb{S}}
\newcommand{\B}{\mathcal{B}}
\newcommand{\Hy}{\mathbb{H}}
\newcommand{\C}{\mathbb{C}}
\newcommand{\psl}{\mbox{\upshape{PSL}}}

\begin{document}

\title[Hyperbolic Structures on the Borromean link complement]{\sc{The exact Computation of a Real Hyperbolic Structure on the Complement of the Borromean Link}  }

\author{Alejandro Ucan-Puc}
\address{ Institut de Mathématiques, Université Pierre et Marie Curie, 4 Place Jussieu, F-75252, Paris, France.}
\email{alejandro.ucan-puc@imj-prg.fr}

  \thanks{Supported by grants of the CONACYT's project 265667} 

%   General info
\subjclass{ }

% \date{January 1, 1994 and, in revised form, June 22, 1994.}

%\dedicatory{This paper is dedicated to .}
\begin{abstract}
In this paper we present a classical construction of the Hyperbolic structure of the complement of a link in the sense of Thurston for the particular case of the Borromean rings link. As this is nothing new, the aim of this paper should be think as to complete the literature about this topic in this particular case.
\end{abstract}

\maketitle 

\section*{Introduction}

The construction of the real hyperbolic structure on the sphere complements of link is a well studied tool introduced by W. Thurston. In Thurston's book (\cite{Thurston1997}) the most detailed construction of this structure is Figure-Eight Knot, but in the literature there are other works that detailed this construction for other knots and links, for example in \cite{Menasco1983}, \cite{Takahashi1985},\cite{Weeks2005}, \cite{Tsvietkova2012}, and \cite{Purcell2020}.

Naively speaking, Thurston's algorithm is reduced to: given a set of tetrahedra in $\overline{\Hy^3}$ and the combinatorial information of how ``glue'' together these tetrahedra, this combinatorial information have to satisfy some regularities to guarantee a ``well-gluing'' (to know a equation system). If we an obtain all these, we can obtain a real hyperbolic structure on the sphere complement of the link, i.e., we will obtain a manifold locally modeled in $\Hy^3$ whose transition charts are elements of $\psl(2,\C).$  The set of tetrahedrons associated to a link complement is induced by the link diagram, the regularities are traduced from the obstructions to obtain a manifold at the moment to glue the tetrahedra. 

It is important to mention that the algorithm complexity highly depends on the number of tetrahedra in the triangulation, with more identified tetrahedra the complexity of the equation system. For this reason, the only examples of hyperbolic structures on the complement of links with exact calculation is made for links with $\leq 4$ tetrahedra in its triangulation. For this reason, we provide a complete example of this construction for the case of the Borromean link that needs a 8 tetrahedrons in its triangulation. 

 Aitchison, Rubinstein and Lumsden (see \cite{Aitchison1992} ) proved that given a general graph there is a pair of simplicial complex such that under identification provide a space that is homeomorphic the closure of the complement of the graph in the sphere; in the case of hyperbolic links this simplicial complex is a perfect triangulation to compute the hyperbolic structure. As other important feature of the present work is to relate the Aitchison-Rubinstein-Lumsden polyhedral decomposition with the generic decomposition used in \cite{Thurston1997}, \cite{Menasco1983} and \cite{Purcell2020}.

The organization of the paper is the following. In section 1 we will explain the Aitchison-Lumsden-Rubinstein polyhedral decomposition for the Borromean links and the gluing data associated to this decomposition. In section 2, we will describe the complex coordinate associated to tetrahedra and the gluing compatibility equations. In Section 3, we will describe the compatibility equations for the particular case of the Borromean link and the solutions of the equation system. And finally, in Section 4, we will compute the representation of the fundamental group of the complement of the Borromean link into $\psl(2,\C).$ 

\section*{Acknowladgements}

I thank Angel Cano for the fruitful discussions over the hyperbolic structures on the complement of hyperbolic links. 
The following paper works as an erratum of the first part of Chapter 5 of \cite{Ucan2019}.

\section{From a Link Diagram to a Manifold}

Let $K$ be a link on $\Es^3$ and $\Gamma_K$ denotes the graph obtained by project $K$ in a plane in $\Es^3.$ Assume that $K$ is such that $\Gamma_K$ is a planar graph, 4-valent and such that the complementary regions admit a checker-board coloring. 

The projection plane decomposes $\Es^3$ into two closed 2-balls (we can think as the upper and lower half spaces in $\mathbb{R}^3$ compactified) and $\Gamma_K$ belongs to the boundary of each of this balls and endows this 3-balls with an abstract polyhedron structure or a CW-complex structure. Let denote by $\Pi^{\pm}_{\Gamma_K}$ to this polyhedra. Notice that we can obtain $\Es^3$ by gluing the corresponding faces of these polyhedra under the identity map. 

From the fact that we assume that the complementary regions of $\Gamma_K$ in the projection plane admit a checker-board coloration, we fixed one for $\Pi^+_{\Gamma_K}$ and give to $\Pi^-_{\Gamma_K}$ the dual coloration. If we change the identification map on the faces of the polyhedra by a  rotation $\sigma(f) \frac{2\pi}{n(f)}$ where $\sigma(f)$ is $+$ if $f$ is white and $-$ otherwise, and $n(f)$ is the number of edges in $f$ for a face $f\in \Pi^+_{\Gamma_K}$ we obtain a space that we will denote by $\overline{X}.$

Lets $X=\overline{X}\setminus \{\mbox{class of the vertices}\}.$

\begin{theorem}[\cite{Aitchison1992}]
\label{ALRTheorem}
The space $X$ is homeomorphic to $\Es^3\setminus K.$
\end{theorem}

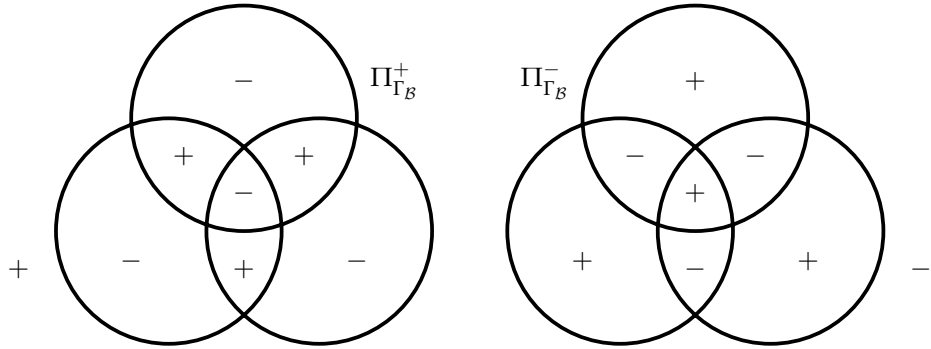
\begin{figure}[H]
\centering
\begin{tikzpicture}
\draw[line width=0.5mm] (-3,1) circle(1.5)
 (-4,-0.5) circle(1.5)
(-2,-0.5) circle(1.5);
\node at (-1,1.5){$\Pi^+_{\Gamma_\B}$};
\node at (-3,1.5){ $-$};
\node at (-3,0){ $-$};
\node at (-3.8,0.5){$+$};
\node at (-2.2,0.5){$+$};
\node at (-3,-1){$+$};
\node at (-4.5,-0.9){$-$};
\node at (-1.5,-0.9){$-$};
\node at (-6,-1){$+$};
\draw[line width=0.5mm] (3,1) circle(1.5)
 (2,-0.5) circle(1.5)
(4,-0.5) circle(1.5);
\node at (1,1.5){$\Pi^-_{\Gamma_\B}$};
\node at (3,1.5){ $+$};
\node at (3,0){ $+$};
\node at (3.8,0.5){$-$};
\node at (2.2,0.5){$-$};
\node at (3,-1){$-$};
\node at (4.5,-0.9){$+$};
\node at (1.5,-0.9){$+$};
\node at (6,-1){$-$};
\end{tikzpicture}
\label{CheckerBoardColoring}
\caption{Plannar projection of $\B$ in $\Es^3.$}
\end{figure}

Lets denote by $\B$ the Borromean rings link and by $\Gamma_\B$ its projection on a plane in $\Es^3.$ The figure \ref{CheckerBoardColoring} provide the coloring of the regions in $\Pi^{\pm}_{\Gamma_\B}.$

It is not hard to picture that the previous polyhedra is similar to an octahedron. Recall that a octahedron, as a simplicial complex, is obtained as the union of four tetrahedra sharing a common edge. At this point, we obtained a triangulation of $\Es^3\setminus \B$ conformed by eight tetrahedrons. 

\section{Ideal Polyhedrons on the Hyperbolic Space}

In the previous section we determine the (combinatorial) polyhedra that become the pieces to glue together and obtain an Euclidean manifold homeomorphic to the complement of a link. Intuitively, if we want an hyperbolic structure on the complement of our link we need to found a model of our polyhedra in the Hyperbolic space $\Hy^3.$ 

\begin{definition}
We will say that a polyhedron $P$ in $\overline{\Hy^3}$ is \emph{ideal} is all its vertices belong to $\partial\Hy^3.$
\end{definition}

Notice that if $P$ is an ideal polyhedron their edges are bi-infinite geodesics.

\subsection{Ideal Hyperbolic Tetrahedron}
The most simple polyhedron in $\Hy^3$ with volume not zero is the \emph{ideal tetrahedron}, see figure \ref{IdealTetrahedron}. 

\tdplotsetmaincoords{60}{110}

%define polar coordinates for some vector
%TODO: look into using 3d spherical coordinate system
\pgfmathsetmacro{\radius}{1}
\pgfmathsetmacro{\thetavec}{0}
\pgfmathsetmacro{\phivec}{0}
\begin{figure}[h]
 \centering
 \begin{tikzpicture}[tdplot_main_coords]
 \tdplotsetthetaplanecoords{\phivec} 
  \draw[tdplot_rotated_coords] (\radius,-1.35,-0.85) arc (-90:70:1.45);
   \tdplotsetrotatedcoords{-5}{-90}{0}

  \draw[tdplot_rotated_coords] (\radius,-3.2,-0.56) arc (-90:91:1.02);
%  \tdplotsetrotatedcoords{-50}{-90}{0}
%  \draw[tdplot_rotated_coords] (\radius -0.1,-0.7,0.19) arc (-90:100:0.5);
%  \draw[tdplot_rotated_coords] (\radius +0.5,-0.7,-1.3  ) arc (-90:100:0.5);
 \tdplotsetrotatedcoords{40}{-90}{5} 
  \draw[tdplot_rotated_coords] (2,-3.5,0) arc (-92:90:1.47);
%  \draw[tdplot_rotated_coords] (\radius+0.9,-1.75,-0.8) arc (-90:100:0.62);
 
   \draw (\radius,0,2.5)--(\radius,0,6);
  \draw (\radius,-1,1.35)--(\radius,-1,5.8);
  \draw (\radius,-3,1.4)--(\radius,-3,5.5);
   \draw[fill=blue] (\radius,-3,1.35) circle(1.8pt) node[below]{$z$};
  \draw[fill=red] (\radius,-1,1.3) circle(1.8pt) node[below]{$1$};
  \draw[fill=green!80!blue] (\radius,0,2.5) circle(1.8pt) node[below]{$0$};
  \end{tikzpicture}
  \label{IdealTetrahedron}
  \caption{Ideal Tetrahedron in $\Hy^3$.}

\end{figure}
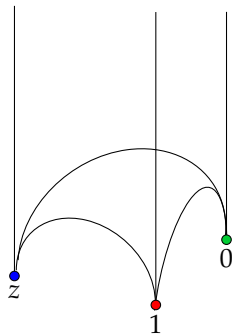
  
  Notice that if we take the intersection of a small horosphere based in one vertex and the ideal tetrahedron, we obtain an Euclidean triangle; even more, if we ask that for the tetrahedron to be oriented, this triangle is well-defined up to orientation preserving similarities. Therefore, an oriented ideal tetrahedron is determined up to orientation preserving similarities by the properties of 4 Euclidean triangles (one in each vertex).

If we take the link triangle for the vertex at infinity in the figure \ref{IdealTetrahedron}, we notice that its angles coincide with the dihedral angles of the polyhedron. So, there are 12 angles (three for each triangle) but as the polyhedron is fixed, two triangles have one angle in common if the vertices of which are associated are in a same edge. Therefore, an ideal Tetrahedron in $\Hy^3$ is determined by three real angles and the angles are paired as follows:

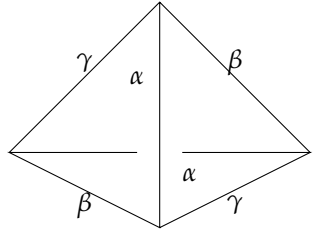
\begin{figure}[h]
\centering
\begin{tikzpicture}
\draw (0,2)--(0,-1)--(-2,0)--(-0.3,0)
(0.3,0)--(2,0)--(0,-1)
(2,0)--(0,2)--(-2,0);
\node at (0.4,-0.3) {$\alpha$};
\node at (-0.3,1){$\alpha$};
\node at (1,1.2){$\beta$};
\node at (-1,1.2){$\gamma$};
\node at (1,-0.7){$\gamma$};
\node at (-1,-0.7){$\beta$};
\end{tikzpicture}
\label{DihedralTetrahedron}
\caption{Pairing of the dihedral angles of a Tetrahedron.}
\end{figure}

Nevertheless if we only depend on the real angles to characterize an hyperbolic ideal tetrahedron, we can found some inconsistencies at the time of gluing. For example Weeks in \cite{Weeks2005} comment that in the gluing process these triangles need to form a $2\pi$ angle and the lengths of the triangles need to coincide (see figure \ref{GluingInconsistences}). In order to prevent this inconsistencies we will change this dihedral real angles by \emph{complex dihedral angles}, for each angle $\alpha$ ($\beta, \,\gamma$) we associate the number $z_0$ ($z_1,\,z_2$) such that $\arg(z_0)=\alpha$ and $\vert z_0\vert$ coincides with the ratio of the side lengths in the triangle. A geometric interpretation of the complex dihedral angles is the complex number that rotates counterclockwise one side of the triangle to an adjacent side. 

\begin{figure}[h]
\begin{subfigure}{.3\textwidth}
\begin{tikzpicture}[scale=0.5]
\draw(0,0)--(3,0)--(1,3)--(-3,2)--(-1,-2)--(0,0);
\node at (1.5,0.2){$z_0$}; 
\draw(0,0)--(1,3) node[pos=0.5,left]{$z_2$};
\draw(0,0)--(-3,2) node[pos=0.4,left]{$z_3$};
\draw(0,0)--(-1,-2) node[pos=0.5,right]{$z_4$};
\end{tikzpicture}
\caption{Angle inconsistency.}
\end{subfigure}
\begin{subfigure}{.3\textwidth}
\begin{tikzpicture}[scale=0.5]
\draw(0,0)--(3,0)--(1,3)--(-3,2)--(-1,-2)--(4.5,0)--(3,0);
\node at (1.5,0.2){$z_0$}; 
\node at (2.7,-0.2){$z_5$};
\draw(0,0)--(1,3) node[pos=0.5,left]{$z_2$};
\draw(0,0)--(-3,2) node[pos=0.4,left]{$z_3$};
\draw(0,0)--(-1,-2) node[pos=0.5,right]{$z_4$};
\end{tikzpicture}
\caption{Modulus inconsistency.}
\end{subfigure}

\begin{subfigure}{.3\textwidth}
\begin{tikzpicture}[scale=0.5]
\draw(0,0)--(3,0)--(1,3)--(-3,2)--(-1,-2)--(3,0);
\node at (1.5,0.2){$z_0$}; 
\node at (1.5,-0.2){$z_5$};
\draw(0,0)--(1,3) node[pos=0.5,left]{$z_2$};
\draw(0,0)--(-3,2) node[pos=0.4,left]{$z_3$};
\draw(0,0)--(-1,-2) node[pos=0.5,right]{$z_4$};
\end{tikzpicture}
\caption{Edge class without inconsistences.}
\end{subfigure}
\label{GluingInconsistencies}
\end{figure}

Let regard the plane containing a triangle vertex as $\mathbb{C}$ and let the triangle have vertices $(v,u,t),$ then the complex dihedral angles are:
\begin{equation}
z(v)=\frac{t-v}{u-v},\quad z(t)=\frac{u-t}{v-t}\quad z(u)=\frac{v-u}{t-u}.
\end{equation}

A direct computation shows that $z(v)z(t)z(u)=-1$ and $z(t)=\frac{1}{1-z(u)}$ and $z(u)=1-\frac{1}{z(u)}.$

\subsection{Triangulations and Gluings}

\begin{definition}
Let $M$ be a 3-manifold. A \emph{topological ideal triangulation} of $M$ is a combinatorial way of gluing ideal tetrahedra so that the result is homeomorphic to $M.$ The ideal vertices correspond to the boundary of $M.$ We also assume that the gluings should take faces to faces, edges to edges, etc.
\end{definition}

Some examples of triangulations are in \cite{Thurston1997} of the Figure-eigth knot, \cite{Menasco1983} of other knots and triangles and the Theorem stated at the first section is a generalization of this definition.

As intuively mention in the previous subsection, we need consistency at the time to glue tetrahedra in the hyperbolic space: that around an edge we complete a $2\pi$ angle and that the sides of the vertex triangles coincide. This can be translated in the following theorem

\begin{theorem}[\cite{Purcell2020}]
Let $M^3$ admit a topological ideal triangulation such that each ideal tetrahedron has a hyperbolic structure. The hyperbolic structures on the ideal tetrahedra induce a hyperbolic structure on the gluing if and only if for each edge class $\mathbf{e}=(e_1,\cdots, e_n)$ we have 
\begin{equation}
\prod_{i=1}^n z(e_i) =1 \qquad \sum_{i=1}^n \arg(z(e_i))=2\pi
\end{equation}
\end{theorem}

If the previous equations system has a solution, we can guarantee the existence of a real hyperbolic structure on $M^3.$ We recall that in the case of geometric structures on manifolds, the more interesting examples are the ones called \emph{complete}, that intituively allow us to obtain $M$ as a quotient of the geometric space by a group isomorphic to the fundamental group of the original manifold.  

\begin{definition}
Let $M$ be a 3-manifold with torus boundary. The \emph{cusp} of $M$ is defined as the neighborhood of $\partial M$ homeomorphic to $T^2\times I$ where $I$ is the interval and $T^2$ is a torus. We will call to each torus component of $\partial M$ a \emph{cusp torus}.
\end{definition}

If $M^3$ has a topological triangulation, we can induce a topological triangulation on each cusp torus boundary by taking the truncated vertices of each ideal tetrahedron. The complete triangulation of the cusps gives a fundamental region of the cusp torus. We are interested in the case where the gluings induce an Euclidean structure on the cusps, \emph{i.e.,} we would like that the cusp torus be flat torus. We have to recall, that the fundamental group of the cusp torus can be obtained by the gluings, the image of this fundamental groups is called \emph{Boundary Holonomy}, and our expectations for the flatness of the cusp torus is translated in that the boundary holonomy is a group generated by Euclidean translations. 

Assume that $\alpha\in\pi_1(T)$ of some cusp torus of $M^3$ and the fundamental domain of the cusp $T$ is such that our first triangle has vertices in $0,\,1\,z(e_1).$ We have that $\alpha$ is an Euclidean traslation if at the moment to change the triangle side from 0 to 1 through all the triangulation, at the end this vector points in the same direction and have the same lenght.

\begin{definition}
Let $[\alpha]\in \pi_1(T)$ a homotopy class and $\alpha$ be representative of the class such that is oriented loop and monotonically, \emph{i.e.,} can be homotoped to a loop that pass through any triangle corner only one time. We will associate the number $H(\alpha)$ as follows:
Let $z_1,\cdots, z_n$ the complex dihedral angles of the corners that cut the loop and to each corner associate $\epsilon_i=\pm 1$ where $\epsilon_i$ is positive if the corner is in the left side of $\alpha$ and negative in the other case. Finally, $H(\alpha)$ is equal to \[\prod_{i=1}^n z_i^{\epsilon_i}\]
\end{definition}

Therefore with this definition, the previous paragraph translates to the equation $H(\alpha)=1$ for every generator of the fundamental group of each cusp.

\section{Borromean Link Triangulation}

As we state in the first section, the complement of the Borromean link can be obtained by the identification of two ideal octahedrons. Nevertheless, the equation system needs a triangulation of the Borromean link's complement, therefore from the Aitchison-Lumsden-Rubinstein polyhedron we will construct a triangulation for the Borromean link complement.

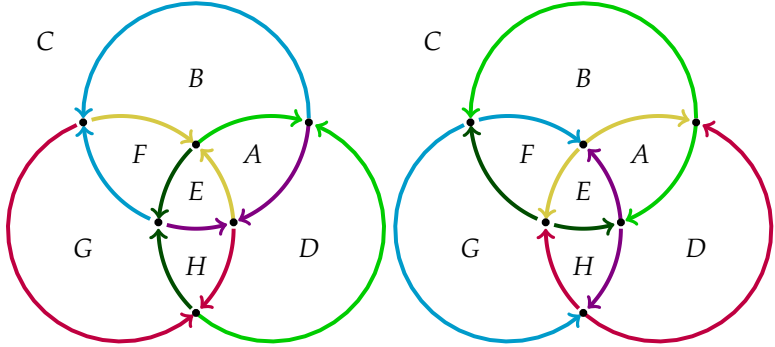
\begin{figure}[h]
 \centering
 \begin{minipage}{0.4\linewidth}
   \begin{tikzpicture}[line width=1.5pt,scale=0.5]
 %Up Circle
 \draw[->,cyan!80!black,domain=-1.1:181.1] plot({3*cos(\x)},{1+3*sin(\x)});
  \draw[->,violet,domain=-5:-68] plot({3*cos(\x)},{1+3*sin(\x)});
 \draw[->,violet,domain=255:286] plot({3*cos(\x)},{1+3*sin(\x)});
  \draw[->,cyan!80!black,domain=-114:-175] plot({3*cos(\x)},{1+3*sin(\x)});
 
 %Right Circle
 \draw[->,green!80!black,domain=-130:67] plot({2+3*cos(\x)},{-2+3*sin(\x)}) ;
 \draw[->,green!80!black,domain=-230:-286] plot({2+3*cos(\x)},{-2+3*sin(\x)});
 \draw[->,green!30!black,domain=135:175] plot({2+3*cos(\x)},{-2+3*sin(\x)});
 \draw[->,green!30!black,domain=-135:-180] plot({2+3*cos(\x)},{-2+3*sin(\x)});
 
 %Left Circle
 \draw[->,purple,domain=113:310] plot({-2+3*cos(\x)},{-2+3*sin(\x)}) ;
 \draw[->,purple,domain=0:-45] plot({-2+3*cos(\x)},{-2+3*sin(\x)});
 \draw[->,yellow!80!black,domain=5:45] plot({-2+3*cos(\x)},{-2+3*sin(\x)});
 \draw[->,yellow!80!black,domain=-255:-310] plot({-2+3*cos(\x)},{-2+3*sin(\x)});
 
 %Points
 \filldraw (1,-1.83) circle(1.5pt)
 (-3,0.83) circle(1.5pt) 
 (-1,-1.83) circle(1.5pt)
 (3,0.83) circle(1.5pt) 
 (0,0.24) circle(1.5pt) 
 (0,-4.24) circle(1.5pt);
 
 % Labels of Faces and Edges
 \node at (0,2){$B$};
 \node at (0,-1){$E$};
 \node at (-1.5,0){$F$};
 \node at (1.5,0){$A$};
 \node at (0,-3){$H$};
 \node at (-3,-2.5){$G$};
 \node at (3,-2.5){$D$};
 \node at (-4,3){$C$};
  \end{tikzpicture}
 \end{minipage}
\begin{minipage}{0.4\linewidth}
  \begin{tikzpicture}[line width=1.5pt,scale=0.5]
 %Up Circle
 \draw[->,green!80!black,domain=-1.1:181.1] plot({3*cos(\x)},{1+3*sin(\x)});
  \draw[->,green!80!black,domain=-5:-68] plot({3*cos(\x)},{1+3*sin(\x)});
 \draw[->,green!30!black,domain=255:286] plot({3*cos(\x)},{1+3*sin(\x)});
  \draw[->,green!30!black,domain=-114:-175] plot({3*cos(\x)},{1+3*sin(\x)});
 
 %Right Circle
 \draw[->,purple,domain=-130:67] plot({2+3*cos(\x)},{-2+3*sin(\x)}) ;
 \draw[->,yellow!80!black,domain=-230:-286] plot({2+3*cos(\x)},{-2+3*sin(\x)});
 \draw[->,yellow!80!black,domain=135:175] plot({2+3*cos(\x)},{-2+3*sin(\x)});
 \draw[->,purple,domain=-135:-180] plot({2+3*cos(\x)},{-2+3*sin(\x)});
 
 %Left Circle
 \draw[->,cyan!80!black,domain=113:310] plot({-2+3*cos(\x)},{-2+3*sin(\x)}) ;
 \draw[->,violet,domain=0:-45] plot({-2+3*cos(\x)},{-2+3*sin(\x)});
 \draw[->,violet,domain=5:45] plot({-2+3*cos(\x)},{-2+3*sin(\x)});
 \draw[->,cyan!80!black,domain=-255:-310] plot({-2+3*cos(\x)},{-2+3*sin(\x)});
 
 %Points
 \filldraw (1,-1.83) circle(1.5pt) 
 (-3,0.83) circle(1.5pt) 
 (-1,-1.83) circle(1.5pt)
 (3,0.83) circle(1.5pt) 
 (0,0.24) circle(1.5pt) 
 (0,-4.24) circle(1.5pt);
 
 % Labels of Faces and Edges
 \node at (0,2){$B$};
 \node at (0,-1){$E$};
 \node at (-1.5,0){$F$};
 \node at (1.5,0){$A$};
 \node at (0,-3){$H$};
 \node at (-3,-2.5){$G$};
 \node at (3,-2.5){$D$};
 \node at (-4,3){$C$};
  \end{tikzpicture}
  \end{minipage}
  \caption{Top and bottom polyhedra and gluing data from Theorem \ref{ALRTheorem}.}
  \label{FaceEdgeGluing}
\end{figure}

As previously mention, the polyhedra obtained from Theorem \ref{ALRTheorem} is an octahedron. We can decompose this octahedron into four tetrahedron that share a common edge that joins the top and bottom of the pyramids. With this we obtain a triangulation of the Borromean link complement. 

\begin{figure}[h]
\centering
\begin{minipage}{.45\textwidth}
\includegraphics[scale=0.2]{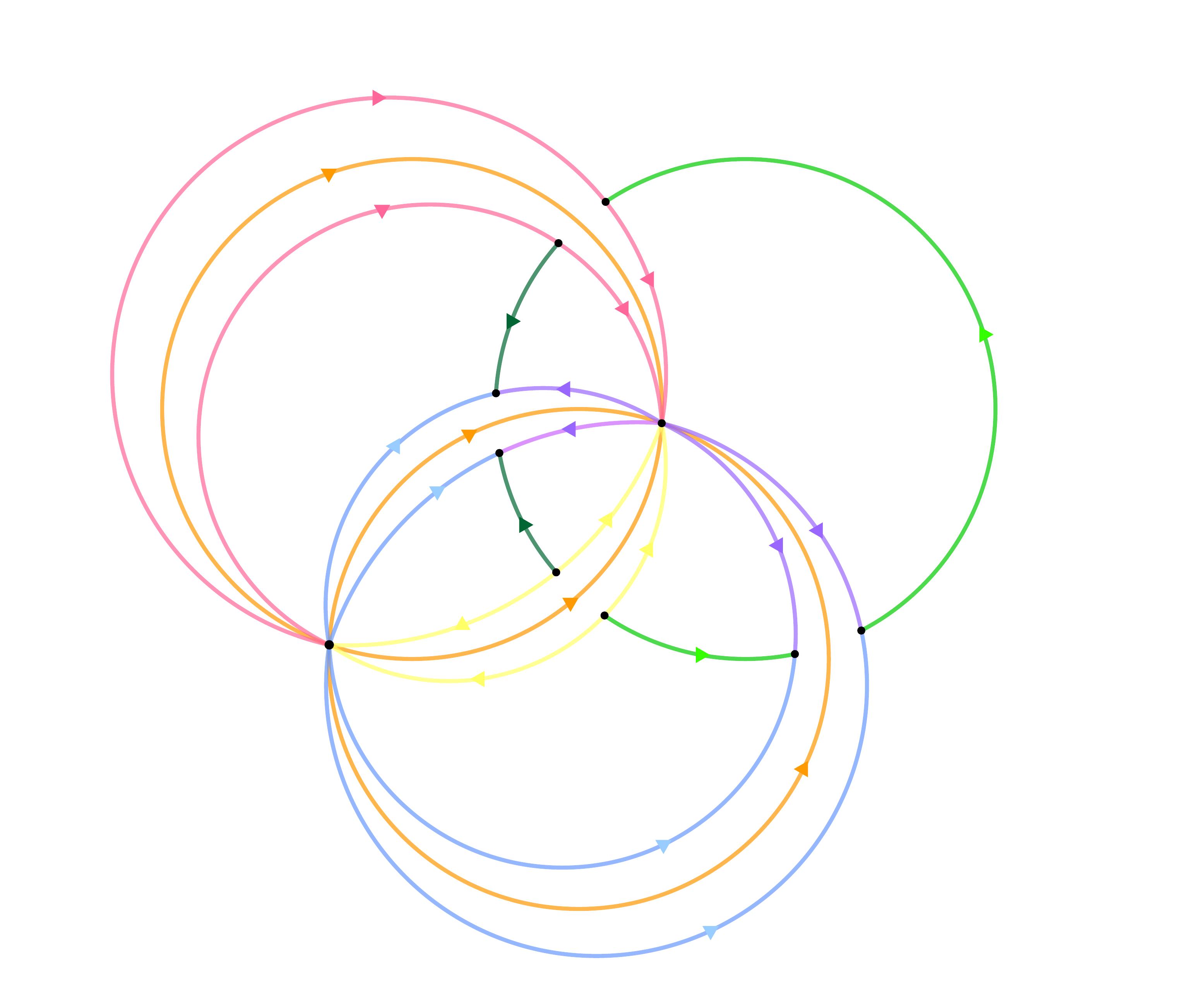}
\end{minipage}
\begin{minipage}{.45\textwidth}
\includegraphics[scale=0.2]{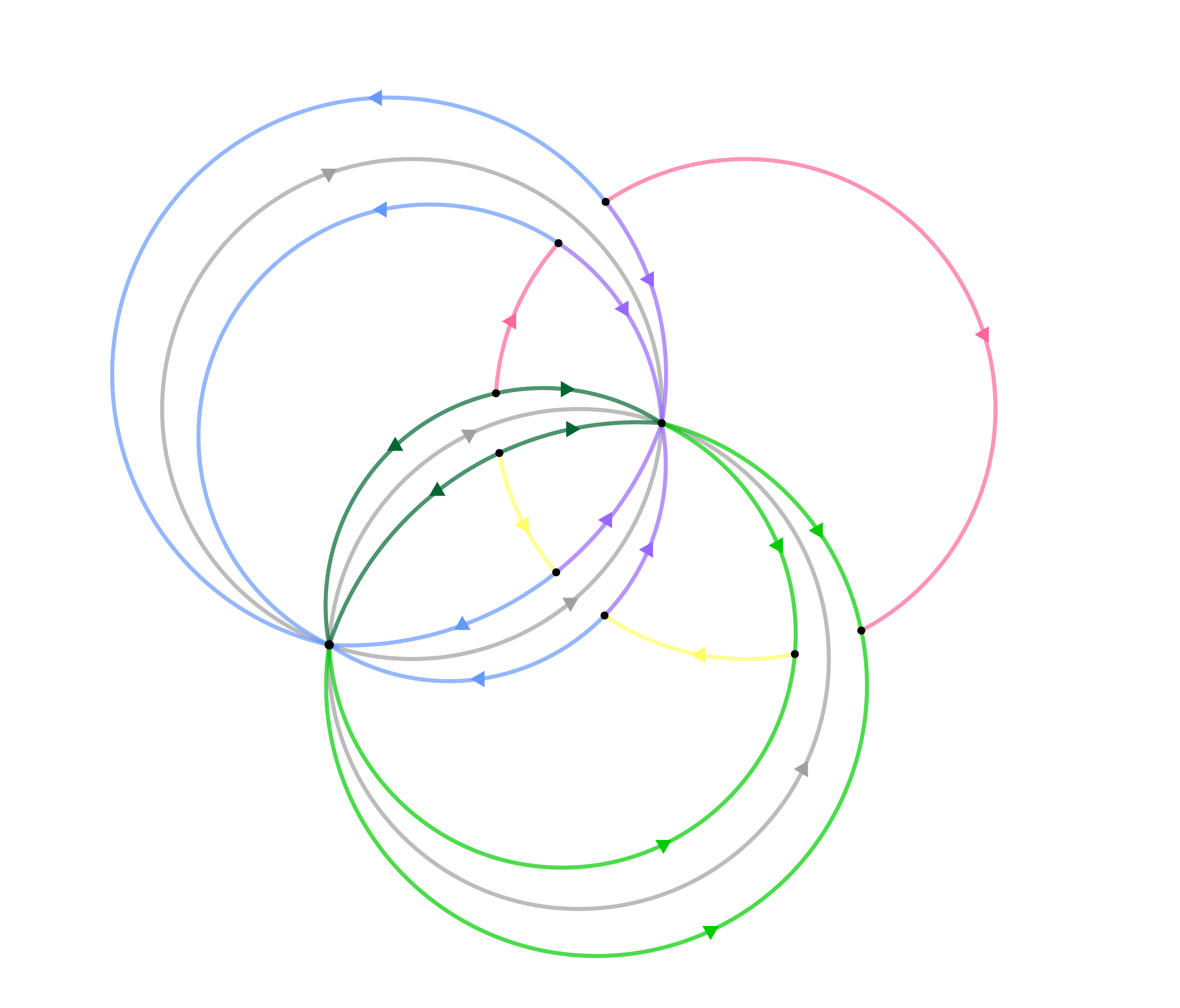}
\end{minipage}
\caption{Top and bottom augmented polyhedron and gluing data.}
\end{figure}

\begin{remark}
The augmented polyhedron are oriented such that coincide with the positively oriented octahedron. Also, from \cite{SnapPy} we know that all tetrahedra in the triangulation is positively oriented.
\end{remark}

The following figure gives in detail all the information of the tetrahedra. To notice, we have 3 classes of vertex and 8 classes of edges.  The vertex-classes correspond to the cusp associated to each component of the link. Also, each tetrahedron has the complex dihedral angles as follow: 
\begin{center}
\begin{tikzpicture}
\draw (0,.1)--(0,1.9) ;
\draw (.1,1.9)--(0.9,1.1);
\draw (0.1,0.1)--(0.9,0.9);
\draw (-0.1,0.1)--(-0.9,0.9);
\draw (-0.1,1.9)--(-0.9,1.1);
\draw (0.1,1)--(0.9,1);
\draw (-0.9,1)--(-0.1,1);
\filldraw[black]  (0,2) circle(1pt) node[above]{\small{1}};
\filldraw[black]  (0,0) circle(1pt) node[below]{\small{0}};
\filldraw[black]  (1,1) circle(1pt) node[right]{\small{3}};
\filldraw[black]  (-1,1) circle(1pt) node[left]{\small{2}};
\node at (0.2,1.2) {$z_i$};
\node at (0.7,1.7) {$v_i$};
\node at (-0.7,1.7){$w_i$};
\node at (-0.7,0.3) {$v_i$};
\node at (0.7,0.3){$w_i$};
\node at (-1.5,2){$T_i$};
\end{tikzpicture}
\end{center}

\begin{figure}[h]
\begin{minipage}{.21\textwidth}
\begin{tikzpicture}
\draw[->, yellow!80!black] (0,.1)--(0,1.9);
\draw[->,cyan!80!black] (.1,1.9)--(0.9,1.1);
\draw[->,green!80!black] (0.1,0.1)--(0.9,0.9);
\draw[->,yellow!80!black] (-0.1,0.1)--(-0.9,0.9);
\draw[->,orange] (-0.1,1.9)--(-0.9,1.1);
\draw[->,violet] (0.1,1)--(0.9,1);
\draw[violet] (-0.9,1)--(-0.1,1);
\filldraw[black]  (0,2) circle(1pt) node[above]{\small{1}};
\filldraw[black]  (0,0) circle(1pt) node[below]{\small{0}};
\filldraw[black]  (1,1) circle(1pt) node[right]{\small{3}};
\filldraw[black]  (-1,1) circle(1pt) node[left]{\small{2}};
\node at (-1.5,2){$T_0$};
\end{tikzpicture}
\end{minipage}
\begin{minipage}{.21\textwidth}
\begin{tikzpicture}
\draw[->, yellow!80!black] (0,.1)--(0,1.9);
\draw[->,orange] (.1,1.9)--(0.9,1.1);
\draw[->,yellow!80!black] (0.1,0.1)--(0.9,0.9);
\draw[->,green!30!black] (-0.1,0.1)--(-0.9,0.9);
\draw[->,cyan!80!black] (-0.1,1.9)--(-0.9,1.1);
\draw[violet] (0.1,1)--(0.9,1);
\draw[<-,violet] (-0.9,1)--(-0.1,1);
\filldraw[black]  (0,2) circle(1pt) node[above]{\small{1}};
\filldraw[black]  (0,0) circle(1pt) node[below]{\small{0}};
\filldraw[black]  (1,1) circle(1pt) node[right]{\small{3}};
\filldraw[black]  (-1,1) circle(1pt) node[left]{\small{2}};
\node at (-1.5,2){$T_1$};
\end{tikzpicture}
\end{minipage}
\begin{minipage}{.21\textwidth}
\begin{tikzpicture}
\draw[<-, cyan!80!black] (0,.1)--(0,1.9);
\draw[->,orange] (.1,1.9)--(0.9,1.1);
\draw[<-,violet] (0.1,0.1)--(0.9,0.9);
\draw[->,green!30!black] (-0.1,0.1)--(-0.9,0.9);
\draw[->,purple] (-0.1,1.9)--(-0.9,1.1);
\draw[purple] (0.1,1)--(0.9,1);
\draw[<-,purple] (-0.9,1)--(-0.1,1);
\filldraw[black]  (0,2) circle(1pt) node[above]{\small{1}};
\filldraw[black]  (0,0) circle(1pt) node[below]{\small{0}};
\filldraw[black]  (1,1) circle(1pt) node[right]{\small{3}};
\filldraw[black]  (-1,1) circle(1pt) node[left]{\small{2}};
\node at (-1.5,2){$T_2$};
\end{tikzpicture}
\end{minipage}
\begin{minipage}{.21\textwidth}
\begin{tikzpicture}
\draw[<-, orange] (0,.1)--(0,1.9);
\draw[->,cyan!80!black] (.1,1.9)--(0.9,1.1);
\draw[->,violet] (0.1,0.1)--(0.9,0.9);
\draw[->,purple] (-0.1,0.1)--(-0.9,0.9);
\draw[->,purple] (-0.1,1.9)--(-0.9,1.1);
\draw[->,green!80!black] (0.1,1)--(0.9,1);
\draw[green!80!black] (-0.9,1)--(-0.1,1);
\filldraw[black]  (0,2) circle(1pt) node[above]{\small{1}};
\filldraw[black]  (0,0) circle(1pt) node[below]{\small{0}};
\filldraw[black]  (1,1) circle(1pt) node[right]{\small{3}};
\filldraw[black]  (-1,1) circle(1pt) node[left]{\small{2}};
\node at (-1.5,2){$T_3$};
\end{tikzpicture}
\end{minipage}\\

\begin{minipage}{.21\textwidth}
\begin{tikzpicture}
\draw[->, green!80!black] (0,.1)--(0,1.9);
\draw[->,gray] (.1,1.9)--(0.9,1.1);
\draw[->,green!80!black] (0.1,0.1)--(0.9,0.9);
\draw[<-,purple] (-0.1,0.1)--(-0.9,0.9);
\draw[<-,cyan!80!black] (-0.1,1.9)--(-0.9,1.1);
\draw[->,violet] (0.1,1)--(0.9,1);
\draw[violet] (-0.9,1)--(-0.1,1);
\filldraw[black]  (0,2) circle(1pt) node[above]{\small{1}};
\filldraw[black]  (0,0) circle(1pt) node[below]{\small{0}};
\filldraw[black]  (1,1) circle(1pt) node[right]{\small{3}};
\filldraw[black]  (-1,1) circle(1pt) node[left]{\small{2}};
\node at (-1.5,2){$T_4$};
\end{tikzpicture}
\end{minipage}
\begin{minipage}{.21\textwidth}
\begin{tikzpicture}
\draw[->, violet] (0,.1)--(0,1.9);
\draw[<-,green!80!black] (.1,1.9)--(0.9,1.1);
\draw[<-,yellow!80!black] (0.1,0.1)--(0.9,0.9);
\draw[->,cyan!80!black] (-0.1,0.1)--(-0.9,0.9);
\draw[<-,gray] (-0.1,1.9)--(-0.9,1.1);
\draw[green!80!black] (0.1,1)--(0.9,1);
\draw[<-,green!80!black] (-0.9,1)--(-0.1,1);
\filldraw[black]  (0,2) circle(1pt) node[above]{\small{1}};
\filldraw[black]  (0,0) circle(1pt) node[below]{\small{0}};
\filldraw[black]  (1,1) circle(1pt) node[right]{\small{3}};
\filldraw[black]  (-1,1) circle(1pt) node[left]{\small{2}};
\node at (-1.5,2){$T_5$};
\end{tikzpicture}
\end{minipage}
\begin{minipage}{.21\textwidth}
\begin{tikzpicture}
\draw[<-, gray] (0,.1)--(0,1.9);
\draw[->,green!30!black] (.1,1.9)--(0.9,1.1);
\draw[->,green!30!black] (0.1,0.1)--(0.9,0.9);
\draw[<-,violet] (-0.1,0.1)--(-0.9,0.9);
\draw[<-,cyan!80!black] (-0.1,1.9)--(-0.9,1.1);
\draw[->,purple] (0.1,1)--(0.9,1);
\draw[purple] (-0.9,1)--(-0.1,1);
\filldraw[black]  (0,2) circle(1pt) node[above]{\small{1}};
\filldraw[black]  (0,0) circle(1pt) node[below]{\small{0}};
\filldraw[black]  (1,1) circle(1pt) node[right]{\small{3}};
\filldraw[black]  (-1,1) circle(1pt) node[left]{\small{2}};
\node at (-1.5,2){$T_6$};
\end{tikzpicture}
\end{minipage}
\begin{minipage}{.21\textwidth}
\begin{tikzpicture}
\draw[->, cyan!80!black] (0,.1)--(0,1.9);
\draw[->,green!30!black] (.1,1.9)--(0.9,1.1);
\draw[<-,yellow!80!black] (0.1,0.1)--(0.9,0.9);
\draw[->,violet] (-0.1,0.1)--(-0.9,0.9);
\draw[->,gray] (-0.1,1.9)--(-0.9,1.1);
\draw[->,green!30!black] (0.1,1)--(0.9,1);
\draw[green!30!black] (-0.9,1)--(-0.1,1);
\filldraw[black]  (0,2) circle(1pt) node[above]{\small{1}};
\filldraw[black]  (0,0) circle(1pt) node[below]{\small{0}};
\filldraw[black]  (1,1) circle(1pt) node[right]{\small{3}};
\filldraw[black]  (-1,1) circle(1pt) node[left]{\small{2}};
\node at (-1.5,2){$T_7$};
\end{tikzpicture}
\end{minipage}

\label{TriangulationBorromean}
\caption{Triangulation of the Borromean's link complement.}
\end{figure}

\subsection{Equations System} \label{SS1:EquationSystem}

\begin{itemize}
\item Edge Equations

\begin{eqnarray*}
z_0 w_1 z_2 w_3 z_4 z_5 v_6 v_7 & = & 1 \\
v_4 w_5 z_6 w_7 & = & 1\\
w_0 v_1 v_2 z_3 & = & 1\\
v_0 z_1 w_2 v_3 w_4 v_5 w_5 z_7 & = & 1\\
v_1 v_2 v_6 w_6 z_7 v_7 & = & 1\\
z_0 z_3 z_4 w_4 z_5 v_5 & = & 1\\
z_0 v_0 z_2 w_2 w_5 w_7 & = &1\\
\end{eqnarray*}
with the restrictions on the sum of their arguments be equal to $2\pi$ for each equation.

\item Boundary Holonomy Equations

\begin{itemize}
\item Vertex Class 0: 0(0), 5(3), 4(0), 3(2), 1(2), 6(3), 7(3), 2(0)

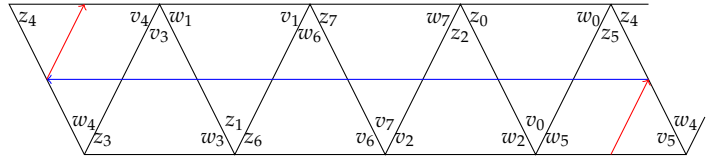
\begin{figure}[h]
\begin{tikzpicture}
\draw (0,0)--(8,0)--(7,2)--(-1,2)--(0,0);
\draw (7,2)--(7.5,2) (8,0)--(8.25,0.5) 
(7,2)--(6,0)--(5,2)--(4,0)--(3,2)--(2,0)--(1,2)--(0,0);
\draw[<-,blue] (-0.5,1)--(7.5,1);
\draw[->,red] (-0.5,1)--(0,2);
\draw[->,red] (7,0)--(7.5,1);
\node[right]  at (-1,1.8){\footnotesize{$z_4$}};
\node[left]  at (1,1.8){\footnotesize{$v_4$}};
\node[below]  at (1,1.8){\footnotesize{$v_3$}};
\node[right]  at (1,1.8){\footnotesize{$w_1$}};
\node[left]  at (3,1.8){\footnotesize{$v_1$}};
\node[below]  at (3,1.8){\footnotesize{$w_6$}};
\node[right]  at (3,1.8){\footnotesize{$z_7$}};
\node[left]  at (5,1.8){\footnotesize{$w_7$}};
\node[below]  at (5,1.8){\footnotesize{$z_2$}};
\node[right]  at (5,1.8){\footnotesize{$z_0$}};
\node[left]  at (7,1.8){\footnotesize{$w_0$}};
\node[below]  at (7,1.8){\footnotesize{$z_5$}};
\node[right]  at (7,1.8){\footnotesize{$z_4$}};
\node[above] at (0,0.2){\footnotesize{$w_4$}};
\node[right] at (0,0.2){\footnotesize{$z_3$}};
\node[left] at (2,0.2){\footnotesize{$w_3$}};
\node[above] at (2,0.2){\footnotesize{$z_1$}};
\node[right] at (2,0.2){\footnotesize{$z_6$}};
\node[left] at (4,0.2){\footnotesize{$v_6$}};
\node[above] at (4,0.2){\footnotesize{$v_7$}};
\node[right] at (4,0.2){\footnotesize{$v_2$}};
\node[left] at (6,0.2){\footnotesize{$w_2$}};
\node[above] at (6,0.2){\footnotesize{$v_0$}};
\node[right] at (6,0.2){\footnotesize{$w_5$}};
\node[left] at (8,0.2){\footnotesize{$v_5$}};
\node[above] at (8,0.2){\footnotesize{$w_4$}};
\end{tikzpicture}
\caption{Triangulation of the cusp torus at vertex class 0.}
\end{figure}
\begin{eqnarray*}
\frac{z_4}{v_5} &=& 1\\
\frac{w_4 w_1 v_7 v_0}{v_3 w_6 z_2 z_5} &=&1
\end{eqnarray*}

\item Vertex Class 1: 0(1), 3(1), 5(0), 2(1), 4(2), 1(1), 0(2), 7(0), 6(2), 3(0), 2(3), 1(3)
\begin{figure}[h]
\begin{tikzpicture}
\draw (0,0)--(4,0)--(8,4)--(4,4)--(0,0) (3,3)--(5,1)
(4,0)--(4,4) (2,0)--(6,4) (6,4)--(6,2)--(2,2)--(2,0) (2,0)--(1,1) (6,4)--(7,3);
\draw[->,red] (4.5,0.5)--(4,0.8)--(3.5,1.5)--(3.4,2)--(3.5,2.5)--(4,3)--(4.5,4);
\draw[->,blue] (6.5,2.5)--(6,2.5)--(5,2)--(4.5,1.5)--(4,1.2)--(3.4,1.4)--(3,2)--(2.45,2.45);
\draw[->,blue] (0.7,0)--(0.4,0.4);
\node[right] at (0.1,0.1){\footnotesize{$z_1$}};
\node[left] at (1.9,0.1){\footnotesize{$w_1$}};
\node[above left] at (2.1,0.2){\footnotesize{$w_3$}};
\node[above right] at (1.9,0.2){\footnotesize{$v_4$}};
\node[right] at (2.1,0.1){\footnotesize{$v_3$}};
\node[left] at (4,0.1){\footnotesize{$z_3$}};
\node[above] at (4.2,0.15){\footnotesize{$w_0$}};
\node[below left] at (1.3,1){\footnotesize{$v_1$}};
\node[right] at (1,1){\footnotesize{$z_3$}};
\node[below ] at (1.9,1.9){\footnotesize{$v_3$}};
\node[right] at (2,1.8){\footnotesize{$w_4$}};
\node[left] at (3.9,1.85){\footnotesize{$z_4$}};
\node[below] at (3.85,1.75){\footnotesize{$w_3$}};
\node[right] at (3.9,1.65){\footnotesize{$z_0$}};
\node[right] at (4.1,1.85){\footnotesize{$z_5$}};
\node[left] at (5.9,1.85){\footnotesize{$v_5$}};
\node[left] at (5,1){\footnotesize{$v_0$}};
\node[above] at (5,1){\footnotesize{$w_5$}};
\node[right] at (2.1,2.1){\footnotesize{$w_6$}};
\node[left] at (3.95,2.15){\footnotesize{$v_6$}};
\node[above] at (3.85,2.15){\footnotesize{$w_1$}};
\node[right] at (3.9,2.35){\footnotesize{$z_2$}};
\node[right] at (4.1,2.1){\footnotesize{$v_7$}};
\node[left] at (6.1,2.15){\footnotesize{$z_7$}};
\node[above] at (6.15,2.1){\footnotesize{$w_2$}};
\node[below] at (3,3){\footnotesize{$z_6$}};
\node[right] at (3,3){\footnotesize{$z_1$}};
\node[left] at (7,3){\footnotesize{$v_2$}};
\node[above] at (7,3){\footnotesize{$w_0$}};
\node[below] at (3.85,3.8){\footnotesize{$v_1$}};
\node[right] at (3.9,3.85){\footnotesize{$v_2$}};
\node[left] at (5.9,3.85){\footnotesize{$w_2$}};
\node[below] at (5.85,3.8){\footnotesize{$w_7$}};
\node[below] at (6.2,3.8){\footnotesize{$z_2$}};
\node[right] at (6.1,3.85){\footnotesize{$z_0$}};
\node[left] at (7.9,3.85){\footnotesize{$v_0$}};
\end{tikzpicture}
\end{figure}

\begin{eqnarray*}
\frac{w_0 v_2}{w_1 v_6 z_4 w_3} & = & 1\\
\frac{z_1 w_6 z_7 w_2}{z_5 z_0 w_3 z_4} & = & 1
\end{eqnarray*}

\item Vertex Class 2: 0(3), 3(3), 5(1), 5(2), 4(1), 4(3), 7(2), 7(1), 6(1), 6(0), 2(2), 1(0)
\begin{figure}[h]
\begin{tikzpicture}
\draw (0,0)--(4,0)--(8,4)--(4,4)--(0,0) (3,3)--(5,1)
(4,0)--(4,4) (2,0)--(6,4) (6,4)--(6,2)--(2,2)--(2,0) (2,0)--(1,1) (6,4)--(7,3);
\draw[->,red] (4.5,0.5)--(4,0.8)--(3.5,1.5)--(3.4,2)--(3.5,2.5)--(4,3)--(4.5,4);
\draw[->,blue] (6.5,2.5)--(6,2.5)--(5,2)--(4.5,1.5)--(4,1.2)--(3.4,1.4)--(3,2)--(2.45,2.45);
\draw[->,blue] (0.7,0)--(0.4,0.4);
\node[right] at (0.1,0.1){\footnotesize{$v_6$}};
\node[left] at (1.9,0.1){\footnotesize{$w_6$}};
\node[above left] at (2.1,0.2){\footnotesize{$z_7$}};
\node[above right] at (1.9,0.2){\footnotesize{$v_2$}};
\node[right] at (2.1,0.1){\footnotesize{$v_7$}};
\node[left] at (4,0.1){\footnotesize{$w_7$}};
\node[above] at (4.2,0.15){\footnotesize{$w_5$}};
\node[below left] at (1.3,1){\footnotesize{$z_6$}};
\node[right] at (1,1){\footnotesize{$w_7$}};
\node[below ] at (1.9,1.9){\footnotesize{$v_7$}};
\node[right] at (2,1.8){\footnotesize{$z_2$}};
\node[left] at (3.9,1.85){\footnotesize{$w_2$}};
\node[below] at (3.85,1.75){\footnotesize{$z_7$}};
\node[right] at (3.9,1.65){\footnotesize{$v_5$}};
\node[right] at (4.1,1.85){\footnotesize{$v_0$}};
\node[left] at (5.9,1.85){\footnotesize{$z_0$}};
\node[left] at (5,1){\footnotesize{$z_5$}};
\node[above] at (5,1){\footnotesize{$w_0$}};
\node[right] at (2.1,2.1){\footnotesize{$w_1$}};
\node[left] at (3.95,2.15){\footnotesize{$z_1$}};
\node[above] at (3.85,2.15){\footnotesize{$w_6$}};
\node[right] at (3.9,2.35){\footnotesize{$w_4$}};
\node[right] at (4.1,2.1){\footnotesize{$v_3$}};
\node[left] at (6.1,2.15){\footnotesize{$w_3$}};
\node[above] at (6.15,2.1){\footnotesize{$z_4$}};
\node[below] at (3,3){\footnotesize{$v_1$}};
\node[right] at (3,3){\footnotesize{$v_6$}};
\node[left] at (7,3){\footnotesize{$v_4$}};
\node[above] at (7,3){\footnotesize{$w_5$}};
\node[below] at (3.85,3.8){\footnotesize{$z_6$}};
\node[right] at (3.9,3.85){\footnotesize{$v_4$}};
\node[left] at (5.9,3.85){\footnotesize{$z_4$}};
\node[below] at (5.85,3.8){\footnotesize{$z_3$}};
\node[below] at (6.2,3.8){\footnotesize{$w_4$}};
\node[right] at (6.1,3.85){\footnotesize{$v_5$}};
\node[left] at (7.9,3.85){\footnotesize{$z_5$}};
\end{tikzpicture}
\end{figure}

\begin{eqnarray*}
\frac{w_5 v_4}{z_7 w_2 z_1 w_6} & = & 1\\
\frac{v_6 z_4 w_3 w_1}{v_0 v_5 z_7 w_2} & = & 1
\end{eqnarray*}
\end{itemize} 
\end{itemize}

\begin{remark}
The cusp triangulations were obtained based in the information given by  software\cite{regina}.
\end{remark}

\subsection{Solutions of the System}
The following complex dihedral angles form a solution of the previous equations system: $z_0=z_5=z_7\frac{1+i}{2},\,z_1=z_2=z_4=1+i,\, z_3=z_6=i.$ We compute this solutions using \cite{SnapPy}, we have to mention that this software algorithm is based in \cite{Weeks2005}.

Using this dihedral angles, we can compute the coordinates of the ideal tetrahedra; even more, using the coordinates of the tetrahedra we can compute the face pairings and the group generated by the face pairings is a representation of the fundamental group of the complement of the Borromean rings into $\psl(2,\C).$

The process to compute the coordinates of a tetrahedron is relatively easy, we need a Lemma stated in \cite{Purcell2020} which relates the computation of the complex dihedral angles after normalization of the tetrahedron, by normalization we mean that the 3 of the vertex of the tetrahedron are $0,\,1,\,\infty.$

\begin{lemma}[Lemma 4.6 in \cite{Purcell2020}]
Let $T$ be an ideal tetrahedron with edge $e_1,$ mapped so that the vertices of $T$ lie at $\infty,0,1$ and $z(e_1).$ The endpoints of $e_1$ lie at $0$ and $\infty.$ Then $T$ has the following addition edge invariants:
\begin{itemize}
\item The edge $e'_1$ opposite to $e_1,$ with vertices $1$ and $z(e_1),$ has edge invariant $z(e'_1)=z(e_1).$
\item The edge $e_2$ with vertices $\infty$ and $1$ has edge invariant $z(e_2)=\frac{1}{1-z(e_1)}.$
\item The edge $e_3$ with vertices $\infty$ and $z(e_1)$ has edge invariant $z(e_3)=\frac{z(e_1)-1}{z(e_1)}.$
\end{itemize}
\end{lemma}

Lets assume that $T_3$ is a normalized tetrahedron, therefore its vertices are $0,\,1,\,\infty$ and $i.$ We will ask that the edge $(01)$ correspond to edge joining $0$ and $\infty,$ therefore the coordinates of the vertices $2$ and $3$ are $i$ and $1$ respectively. Let $\zeta_0(i)$ be the coordinates of $T_0$ for $i=0,1,2,3.$ We know that $T_3$ and $T_1$ are identified by the faces $3(120)$ and $1(123),$ which means that $\zeta_1(1)=\infty,\, \zeta_1(3)=0$ and $\zeta_1(2)=i.$ Using the map $z\mapsto \frac{z}{i},$ me can normalized the tetrahedron $T_1,$ therefore $\frac{\zeta_1(0)}{i}=z(T_0)=1+i$ which implies that $\zeta_1(0)=-1+i.$ Repeating this process for each tetrahedra we can compute all the vertices. Table \ref{TetrahedronVertices} contains an example of vertices obtained by this process.
\begin{table}
\begin{tabular}{|c|c|c|c|c|}
\hline
Tetrahedron & $\zeta_i(0)$ & $\zeta_i(1)$ & $\zeta_i(2)$ & $\zeta_i(3)$ \\ \hline
3& $0$ & $\infty$ & $i$ & $1$ \\ \hline
1 & $-1$ & $\infty$ & $i$ & $0$ \\ \hline
2 & $-i$ & $\infty$ &  $-1$ & $0$ \\ \hline
0 & $-i$ & $\infty$ & $0$ & $1$ \\ \hline
4 & $i$ & $1$ & $\infty$ & $1+2i$ \\ \hline
6 & $1+2i$ & $1$ & $\infty$ & $2+i$ \\ \hline
7 & $1+i$ & $1$ & $1+2i$ & $2+i$ \\\hline
5 & $1+i$ & $1+2i$ & $1$ & $i$ \\ \hline
\end{tabular}
\caption{Coordinates of Tetrahedron in $\partial\Hy^3.$}
\label{TetrahedronVertices}
\end{table}

\section{Fundamental Group Representation.}

Over this section we will compute the representation of the fundamental group associated to the hyperbolic structure computed in the previous section. The fundamental idea is that the group generated by the face-pairing transformations is a representation of the fundamental group, this is due to the well known Poincaré Polyhedron Theorem.

\subsection{Fundamental group of $\Es^3\setminus\B$}

Let $\mathcal{L}$ be a link in $\Es^3,$ and let $\pi_1(\Es^3\setminus\mathcal{L})$ be its fundamental group. There exists an easy way to compute the presentation of the fundamental group on a link complement called the Wirtinger's presentation. 

Intuitively, the Wirtinger's presentation is obtained in the following way: for each directed arc of the link (as a graph with intersections) we take a generator of the fundamental group, and for every crossing we will have a relation between the generators, this relation depends on the type of crossing.

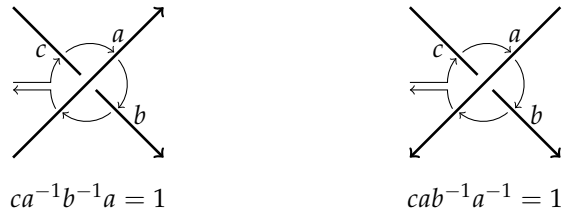
\begin{figure}[h]

\begin{minipage}{.41\textwidth}
\centering
\begin{tikzpicture}
\draw[->,line width=0.35mm] (0,0)--(2,2) node[pos=0.7,above]{$a$};
\draw[->,line width=0.35mm](1.1,.9)--(2,0) node[pos=0.65,above]{$b$};
\draw[line width=0.35mm] (0,2)--(.9,1.1) node[pos=0.65,left]{$c$};
\draw[->] (0,1)--(0.5,1) arc(180:140:0.5);
\draw[->] (0.7,1.4) arc(130:50 :0.5);
\draw[->] (1.4,1.3) arc(40:-40:0.5);
\draw[->] (1.3,0.6) arc(-50:-130:0.5);
\draw[<-] (0,.9)--(0.5,0.9) arc(180:210:0.5);
\node at (1,-0.5){$ca^{-1}b^{-1}a=1$};
\end{tikzpicture}
\end{minipage}
\begin{minipage}{.41\textwidth}
\centering
\begin{tikzpicture}
\draw[<-,line width=0.35mm] (0,0)--(2,2) node[pos=0.7,above]{$a$};
\draw[->,line width=0.35mm](1.1,.9)--(2,0) node[pos=0.65,above]{$b$};
\draw[line width=0.35mm] (0,2)--(.9,1.1) node[pos=0.65,left]{$c$};
\draw[->] (0,1)--(0.5,1) arc(180:140:0.5);
\draw[->] (0.7,1.4) arc(130:50 :0.5);
\draw[->] (1.4,1.3) arc(40:-40:0.5);
\draw[->] (1.3,0.6) arc(-50:-130:0.5);
\draw[<-] (0,.9)--(0.5,0.9) arc(180:210:0.5);
\node at (1,-0.5){$cab^{-1}a^{-1}=1$};
\end{tikzpicture}
\end{minipage}
\label{Wirtingerpresentation}
\caption{Wirtinger's Relations for the different type of crossing.}
\end{figure}

In the particular case of the Borromean link, we have 6 generators and 6 relations, the generators are described in Figure

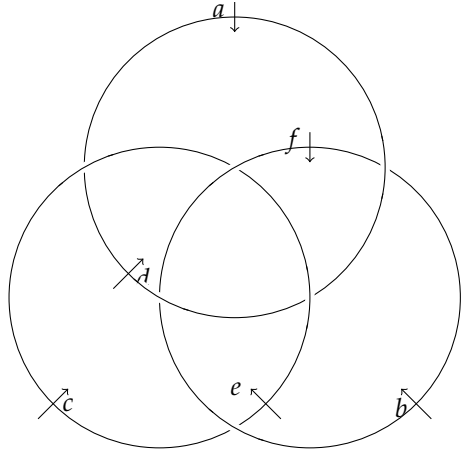
\begin{figure}[h]
\centering
\begin{tikzpicture}
\begin{knot}[
%  draft mode=crossings,
  flip crossing/.list={3,4}
]
\strand (1,0) circle[radius=2cm];
\strand (-1,0) circle[radius=2cm];
\strand (0,{sqrt(3)}) circle[radius=2cm];
\draw[->] (0,{sqrt(3)+2.2})--(0,{sqrt(3)+1.8}) node[pos=0.3,left]{$a$};
\draw[->] (1,2.2)--(1,1.8) node[pos=0.3,left]{$f$};
\draw[->] (-{sqrt(2)-0.2},{sqrt(3)-sqrt(2)-0.2})--(-{sqrt(2)+0.2},{sqrt(3)-sqrt(2)+0.2}) node[below]{$d$};
\draw[->] ({-1-sqrt(2)-0.2},{-sqrt(2)-0.2})--({-1-sqrt(2)+0.2},{-sqrt(2)+0.2}) node[below]{$c$};
\draw[->] ({1+sqrt(2)+0.2},{-sqrt(2)-0.2})--({1+sqrt(2)-0.2},{-sqrt(2)+0.2}) node[below]{$b$};
\draw[->] ({-1+sqrt(2)+0.2},{-sqrt(2)-0.2})--({-1+sqrt(2)-0.2},{-sqrt(2)+0.2}) node[left]{$e$};
\end{knot}
\end{tikzpicture}
\label{WirtingerBorromean}
\caption{Wirtinger Generator and Relations for the Borromean link.}
\end{figure}

Using the generators described in Figure \ref{WirtingerBorromean}, the relations associated to the crossings are:

\begin{eqnarray*}
R_1&=& d^{-1}b a b^{-1},\\
R_2&=& fdc^{-1}d^{-1},\\
R_3&=&efb^{-1}f^{-1},\\
R_4&=&fac^{-1}f^{-1},\\
R_5&=&ecb^{-1}c^{-1},\\
R_6&=&dea^{-1}e^{-1}.
\end{eqnarray*}

From the previous 6 relations, we can deduce that: $d=bab^{-1},$ $f=aca^{-1},$ and $e=cbc^{-1}.$ After this, the relation can be reduced to:

\begin{eqnarray*}
R_2&=&aca^{-1}bab^{-1}c^{-1}ba^{-1}b^{-1},\\
R_3&=&cbc^{-1}aca^{-1}b^{-1}ac^{-1}a^{-1},\\
R_6&=&bab^{-1}cbc^{-1}a^{-1}cb^{-1}c^{-1}.
\end{eqnarray*}

An easy computation proves that $R_6=(R_3 R_2)^{-1}.$ Therefore, from the Wirtinger's relation, the presentation of the fundamental group of the complement of the Borromean link is

\begin{equation}
\pi_1(\Es^3\setminus \B)=\left\langle a,b,c| [c^{-1},[b^{-1},a]]=[b,[c,a^{-1}]]=1\right\rangle.
\end{equation}

\subsection{Face-pairings}

As we said in the introduction of this section, the Poincaré polyhedron theorem allow us to said that the group in $\psl(2,\C)$ of face-pairing transformation is a representation of $\pi_1(\Es^3,\B).$

\begin{theorem}[Poincaré Polyhedron Theorem, \cite{Maskit1988}]
Let $X:=P/\sim$ the quotient space under the equivalence $x\sim \gamma_f(x)$ for $x\in f.$ If $X$ is complete with the quotient metric from $P,$ the map $F:P\to X$ is finite-to-one, and  if for every edge class $\overline{e}\in P_1/\sim$ there exists an integer $N_{\overline{e}}$ such that the sum of the measure dihedral angles incidents to $\overline{e}$ is $\frac{2\pi}{N_{\overline{e}}}.$ Then the group $G\subset \psl(2,\C)$ generated by the set $\left\{\gamma_f:f\in P_2\right\}$ is discrete, $P$ is a fundamental domain of $G$ and $\Hy^3/G$ is homeomorphic to $X.$

\end{theorem}

The algorithm to compute the face-pairing is as follows:
\begin{enumerate}
\item Let $F_1=(i,j,k), F_2=(l,m,n)$ be to faces in Table \ref{TetrahedronVertices} that are paired (see Figure).
\item Let $g_1,g_2\in\psl(2,\C)$ the matrices that send the face $(i,j,k),$ and $(l,m,n)$ (respectively) to $(\infty,0,1).$
\item The face-pairing matrix of $F_1$ and $F_2$ is the matrix $g_2^{-1}g_1.$
 \end{enumerate}

The list of face-pairing transformations are:
\begin{itemize}
\item $3(0\,2\,3)\to 4(2\,0\,3)$ \[\gamma_{34}=\begin{bmatrix}
2i & 1 \\ 1 & 0
\end{bmatrix},\]

\item $0(0\,1\,3)\to 5(3\,0\,2)$ \[\gamma_{05}=\begin{bmatrix}
1+i & -1 \\ 1 & -1+i
\end{bmatrix},\]

\item $0(0\,2\,3)\to 5(3\,0\,1)$ \[\tau_{05}=\begin{bmatrix}
1-2i & 1+i \\ -1-i & 1
\end{bmatrix},\]

\item $1(0\,1\,2)\to 6(1\,2\,3)$ \[\gamma_{16}=\begin{bmatrix}
1 & 2 \\ 0 & 1
\end{bmatrix},\]

\item $1(0\,2\,3)\to 6(0\,3\,2)$ \[\tau_{16}=\begin{bmatrix}
2+2i & 1 \\ 1 & 0
\end{bmatrix},\]

\item $2(0\,1\,2)\to 7(3\,0\,1)$ \[\gamma_{27}=\begin{bmatrix}
2 & 3+i \\ 1-i & 2
\end{bmatrix},\]

\item $2(0\,2\,3)\to 7(3\,2\,0)$ \[\tau_{27}=\begin{bmatrix}
3-3i & 2 \\ -2i & 1-i
\end{bmatrix},\]

\end{itemize}

The generators of the fundamental group are the meridians of the Boundary holonomy group, computed in the section \ref{SS1:EquationSystem}. We have to recall that the curves described in the triangulation are just representatives of the class, we can chose another representative:

\begin{itemize}
\item Meridian in {\textbf{Cusp 0}}  with associated matrix $\gamma_{34}.$
\begin{figure}[h]
\begin{tikzpicture}
\draw (0,0)--(8,0)--(7,2)--(-1,2)--(0,0);
\draw (7,2)--(7.5,2) (8,0)--(8.25,0.5) 
(7,2)--(6,0)--(5,2)--(4,0)--(3,2)--(2,0)--(1,2)--(0,0);
\draw[->,blue] (1,0)--(0,2);
\node[right]  at (-1,1.8){\footnotesize{$z_4$}};
\node[left]  at (1,1.8){\footnotesize{$v_4$}};
\node[below]  at (1,1.8){\footnotesize{$v_3$}};
\node[right]  at (1,1.8){\footnotesize{$w_1$}};
\node[left]  at (3,1.8){\footnotesize{$v_1$}};
\node[below]  at (3,1.8){\footnotesize{$w_6$}};
\node[right]  at (3,1.8){\footnotesize{$z_7$}};
\node[left]  at (5,1.8){\footnotesize{$w_7$}};
\node[below]  at (5,1.8){\footnotesize{$z_2$}};
\node[right]  at (5,1.8){\footnotesize{$z_0$}};
\node[left]  at (7,1.8){\footnotesize{$w_0$}};
\node[below]  at (7,1.8){\footnotesize{$z_5$}};
\node[right]  at (7,1.8){\footnotesize{$z_4$}};
\node[above] at (0,0.2){\footnotesize{$w_4$}};
\node[right] at (0,0.2){\footnotesize{$z_3$}};
\node[left] at (2,0.2){\footnotesize{$w_3$}};
\node[above] at (2,0.2){\footnotesize{$z_1$}};
\node[right] at (2,0.2){\footnotesize{$z_6$}};
\node[left] at (4,0.2){\footnotesize{$v_6$}};
\node[above] at (4,0.2){\footnotesize{$v_7$}};
\node[right] at (4,0.2){\footnotesize{$v_2$}};
\node[left] at (6,0.2){\footnotesize{$w_2$}};
\node[above] at (6,0.2){\footnotesize{$v_0$}};
\node[right] at (6,0.2){\footnotesize{$w_5$}};
\node[left] at (8,0.2){\footnotesize{$v_5$}};
\node[above] at (8,0.2){\footnotesize{$w_4$}};
\end{tikzpicture}
\caption{Representative of a meridian in the cusp 0.}
\end{figure}

\item Meridian in {\textbf{Cusp 0}}  with associated matrix $\gamma_{16}^{-1}.$
\begin{figure}[h]
\begin{tikzpicture}
\draw (0,0)--(4,0)--(8,4)--(4,4)--(0,0) (3,3)--(5,1)
(4,0)--(4,4) (2,0)--(6,4) (6,4)--(6,2)--(2,2)--(2,0) (2,0)--(1,1) (6,4)--(7,3);
%\draw[->,red] (4.5,0.5)--(4,0.8)--(3.5,1.5)--(3.4,2)--(3.5,2.5)--(4,3)--(4.5,4);
\draw[->,red] (1,0)--(2,1)--(3,2)--(2.5,2.5);

\node[right] at (0.1,0.1){\footnotesize{$z_1$}};
\node[left] at (1.9,0.1){\footnotesize{$w_1$}};
\node[above left] at (2.1,0.2){\footnotesize{$w_3$}};
\node[above right] at (1.9,0.2){\footnotesize{$v_4$}};
\node[right] at (2.1,0.1){\footnotesize{$v_3$}};
\node[left] at (4,0.1){\footnotesize{$z_3$}};
\node[above] at (4.2,0.15){\footnotesize{$w_0$}};
\node[below left] at (1.3,1){\footnotesize{$v_1$}};
\node[right] at (1,1){\footnotesize{$z_3$}};
\node[below ] at (1.9,1.9){\footnotesize{$v_3$}};
\node[right] at (2,1.8){\footnotesize{$w_4$}};
\node[left] at (3.9,1.85){\footnotesize{$z_4$}};
\node[below] at (3.85,1.75){\footnotesize{$w_3$}};
\node[right] at (3.9,1.65){\footnotesize{$z_0$}};
\node[right] at (4.1,1.85){\footnotesize{$z_5$}};
\node[left] at (5.9,1.85){\footnotesize{$v_5$}};
\node[left] at (5,1){\footnotesize{$v_0$}};
\node[above] at (5,1){\footnotesize{$w_5$}};
\node[right] at (2.1,2.1){\footnotesize{$w_6$}};
\node[left] at (3.95,2.15){\footnotesize{$v_6$}};
\node[above] at (3.85,2.15){\footnotesize{$w_1$}};
\node[right] at (3.9,2.35){\footnotesize{$z_2$}};
\node[right] at (4.1,2.1){\footnotesize{$v_7$}};
\node[left] at (6.1,2.15){\footnotesize{$z_7$}};
\node[above] at (6.15,2.1){\footnotesize{$w_2$}};
\node[below] at (3,3){\footnotesize{$z_6$}};
\node[right] at (3,3){\footnotesize{$z_1$}};
\node[left] at (7,3){\footnotesize{$v_2$}};
\node[above] at (7,3){\footnotesize{$w_0$}};
\node[below] at (3.85,3.8){\footnotesize{$v_1$}};
\node[right] at (3.9,3.85){\footnotesize{$v_2$}};
\node[left] at (5.9,3.85){\footnotesize{$w_2$}};
\node[below] at (5.85,3.8){\footnotesize{$w_7$}};
\node[below] at (6.2,3.8){\footnotesize{$z_2$}};
\node[right] at (6.1,3.85){\footnotesize{$z_0$}};
\node[left] at (7.9,3.85){\footnotesize{$v_0$}};
\end{tikzpicture}
\caption{Representative of a meridian in the cusp 1.}
\end{figure}

\item Meridian in {\textbf{Cusp 2}}  with associated matrix $\gamma_{16}\gamma_{27}^{-1}(=\gamma_{05}^{-1}).$

\begin{figure}[h]
\begin{tikzpicture}
\draw (0,0)--(4,0)--(8,4)--(4,4)--(0,0) (3,3)--(5,1)
(4,0)--(4,4) (2,0)--(6,4) (6,4)--(6,2)--(2,2)--(2,0) (2,0)--(1,1) (6,4)--(7,3);
%\draw[->,red] (4.5,0.5)--(4,0.8)--(3.5,1.5)--(3.4,2)--(3.5,2.5)--(4,3)--(4.5,4);
\draw[->,red] (3.2,0)--(2.8,0.8)--(3,2)--(3.5,2.5)--(3.5,3.5);

\node[right] at (0.1,0.1){\footnotesize{$z_1$}};
\node[left] at (1.9,0.1){\footnotesize{$w_1$}};
\node[above left] at (2.1,0.2){\footnotesize{$w_3$}};
\node[above right] at (1.9,0.2){\footnotesize{$v_4$}};
\node[right] at (2.1,0.1){\footnotesize{$v_3$}};
\node[left] at (4,0.1){\footnotesize{$z_3$}};
\node[above] at (4.2,0.15){\footnotesize{$w_0$}};
\node[below left] at (1.3,1){\footnotesize{$v_1$}};
\node[right] at (1,1){\footnotesize{$z_3$}};
\node[below ] at (1.9,1.9){\footnotesize{$v_3$}};
\node[right] at (2,1.8){\footnotesize{$w_4$}};
\node[left] at (3.9,1.85){\footnotesize{$z_4$}};
\node[below] at (3.85,1.75){\footnotesize{$w_3$}};
\node[right] at (3.9,1.65){\footnotesize{$z_0$}};
\node[right] at (4.1,1.85){\footnotesize{$z_5$}};
\node[left] at (5.9,1.85){\footnotesize{$v_5$}};
\node[left] at (5,1){\footnotesize{$v_0$}};
\node[above] at (5,1){\footnotesize{$w_5$}};
\node[right] at (2.1,2.1){\footnotesize{$w_6$}};
\node[left] at (3.95,2.15){\footnotesize{$v_6$}};
\node[above] at (3.85,2.15){\footnotesize{$w_1$}};
\node[right] at (3.9,2.35){\footnotesize{$z_2$}};
\node[right] at (4.1,2.1){\footnotesize{$v_7$}};
\node[left] at (6.1,2.15){\footnotesize{$z_7$}};
\node[above] at (6.15,2.1){\footnotesize{$w_2$}};
\node[below] at (3,3){\footnotesize{$z_6$}};
\node[right] at (3,3){\footnotesize{$z_1$}};
\node[left] at (7,3){\footnotesize{$v_2$}};
\node[above] at (7,3){\footnotesize{$w_0$}};
\node[below] at (3.85,3.8){\footnotesize{$v_1$}};
\node[right] at (3.9,3.85){\footnotesize{$v_2$}};
\node[left] at (5.9,3.85){\footnotesize{$w_2$}};
\node[below] at (5.85,3.8){\footnotesize{$w_7$}};
\node[below] at (6.2,3.8){\footnotesize{$z_2$}};
\node[right] at (6.1,3.85){\footnotesize{$z_0$}};
\node[left] at (7.9,3.85){\footnotesize{$v_0$}};
\end{tikzpicture}
\caption{Representative of a meridian in the cusp 2.}
\end{figure}
\end{itemize}

Therefore the elements that generate image of $\pi_1(\Es^3\setminus \B)$ in $\psl(2,\C)$ are: $\gamma_{34},\gamma_{16}^{-1},$ and $\gamma_{05}^{-1}.$ From a direct computation we have:

\begin{eqnarray}
\left[\gamma_{16},\left[\gamma_{34}^{-1},\gamma_{05}\right]\right]&=&1\\
\left[\gamma_{05}^{-1},\left[\gamma_{34},\gamma_{16}^{-1}\right]\right]&=&1
\end{eqnarray}

\bibliographystyle{abbrv}
\bibliography{HRBOR.bib}

\end{document}